\documentclass[10pt,reqno]{amsart}
\usepackage{amssymb}
\usepackage{amsfonts}
\usepackage{mathrsfs}
\usepackage{amsmath}
\usepackage{latexsym}
\usepackage{soul}
\usepackage{float}
\usepackage{tikz}
\usepackage{graphicx}
\usepackage{color}
\usepackage[utf8]{inputenc}
\newtheorem{theorem}{Theorem}
\newtheorem{lemma}{Lemma}

\newcommand{\nl}{\newline}

\newcommand{\R}{{\mathbb{R}}}

\newcommand{\cA}{{\mathcal A}}

\newcommand{\ia}{({\rm i})}
\newcommand{\ib}{({\rm ii})}
\newcommand{\ic}{({\rm iii})}

\newcommand{\be}{\begin{equation}}
\newcommand{\ee}{\end{equation}}
\newcommand{\bea}{\begin{eqnarray}}
\newcommand{\eea}{\end{eqnarray}}
\newcommand{\bean}{\begin{eqnarray*}}
\newcommand{\eean}{\end{eqnarray*}}
\newcommand{\la}{\label}

\newcommand{\cic}{{C^{\infty}_c}}
\newcommand{\sph}{{ {\rm S}^{n-1} }}

\def\Int{\displaystyle\int}

\newcommand{\cico}{C^{\infty}_c(\Omega)}
\newcommand{\mI}{{\mathbb{I}}}
\newcommand{\ino}{\int_{\Omega}}

\begin{document}

\title{Sobolev improvements on sharp Rellich inequalities}

\author{G. Barbatis}
\address{Department of Mathematics, National and Kapodistrian University of Athens, 15784 Athens, Greece}
\email{gbarbatis@math.uoa.gr}

\author{A. Tertikas}
\address{Department of Mathematics and Applied Mathematics, University of Crete, 70013 Heraklion, Crete, Greece \&}
\address{Institute of Applied and Computational Mathematics, FORTH,
71110 Heraklion, Crete, Greece}
\email{tertikas@uoc.gr}


.


\dedicatory{Dedicated to E.B. Davies on the occasion of his 80th birthday}

\keywords{Rellich inequality; Sobolev inequality;
best constant \\
{\em \hspace{2cm} To appear in the Journal of Spectral Theory}}

\begin{abstract}
There are two Rellich inequalities for the bilaplacian,
that is for $\int (\Delta u)^2dx$, the one involving $|\nabla u|$ and the other involving $|u|$ at the RHS.
In this article we consider these inequalities with sharp constants and obtain
sharp Sobolev-type improvements.
More precisely, in our first result we improve 
the Rellich inequality with
$|\nabla u|$ obtained  by Beckner
in dimensions $n=3,4$ by a sharp Sobolev term thus complementing
existing results for the case $n\geq 5$.
In the second theorem the sharp constant of the Sobolev improvement
for the Rellich inequality with $|u|$ is obtained.
\end{abstract}

\maketitle

\section*{Introduction}

The study of PDEs involving the bilaplacian is often related
to functional inequalities for the associated energy, namely
$\int (\Delta u)^2dx$. Two important such inequalities are
the Sobolev inequality and the Rellich inequality.

There are two Rellich inequalities related to the bilaplacian. 
The first one asserts that for $n\geq 5$ there holds
\be
\int_{\R^n}(\Delta u)^2 dx \geq \frac{n^2(n-4)^2}{16}
\int_{\R^n}\frac{u^2}{|x|^4} dx \, , \quad u\in\cic(\R^n) ,
\la{ri1}
\ee
and the constant is the best possible.
Inequality \eqref{ri1} was proved by F. Rellich, see \cite{Rel}. For more results
on inequalities of this type and related improvements
we refer to
\cite{AdS,AGS,B,BT06,CaMu,DMY,DH,DLT,GGM,GeLi,GMP,GM,ST,TZ} and references therein.

The second Rellich inequality is valid not only for $n\geq 5$ but also
for $n=3,4$ and reads
\be
\int_{\R^n}(\Delta u)^2 dx \geq c_n
\int_{\R^n}\frac{|\nabla u|^2}{|x|^2} dx \, , \quad u\in\cic(\R^n) ,
\la{ri2}
\ee
where
\be
c_n=\left\{
\begin{array}{ll}
 \frac{25}{36}, & \quad n=3, \\[0.1cm]
3,  & \quad n=4, \\[0.1cm]
 \frac{n^2}{4}, & \quad  n\geq 5.
\end{array}
\right.
\la{c_n}
\ee
is the best possible constant. Inequality \eqref{ri2} was proved in \cite{TZ}
in case $n\geq 5$ and then by Beckner for any $n\geq 3$ \cite{Be}.
An alternative proof for $n\geq 3$ was given by Cazacu \cite{Ca}.
We note that in cases $n=3,4$ there is a breaking of symmetry. For 
more information on Rellich inequalities in the spirit of 
\eqref{ri2} we refer to \cite{Ca,DMY,DT,KO,TZ}.

The Sobolev inequality for the bilaplacian 
in $\R^n$, $n\geq 5$, reads
\be
\la{1}
\int_{\R^n} (\Delta u)^2 dx \geq S_{2,n} \Big(  
\int_{\R^n} |u|^{\frac{2n}{n-4}}dx\Big)^{\frac{n-4}{n}}  , 
\quad u\in \cic(\R^n).
\ee
The best constant $S_{2,n}$ in \eqref{1} has been
computed in \cite{EFJ} and is given by
\[
S_{2,n} = \pi^2 (n^2-4n)(n^2-4)\Big( \frac{\Gamma(\frac{n}{2})}{\Gamma(n)} \Big)^4.
\]

The aim of this work is to improve the above Rellich inequalities by adding a Sobolev-type term.
In \cite{TZ} improved versions of \eqref{ri1} and
\eqref{ri2} were obtained for a bounded domain $\Omega\subset\R^n$, $n\geq 5$.
More precisely, let $X(r)=(1-\log r)^{-1}$, $0<r<1$, and $D=\sup_{\Omega}|x|$. 
In \cite[Theorem 1.1]{TZ} it was shown that for $n\geq 5$ there exist constants $C_n$ and
$C_n'$ which depend only on $n$ such that for any $u\in\cic(\Omega)$ there holds
\be
\int_{\Omega}(\Delta u)^2 dx - \frac{n^2(n-4)^2}{16}
\int_{\Omega}\frac{u^2}{|x|^4} dx 
\geq C_n \Big(  
\int_{\Omega} X(|x|/D)^{\frac{2(n-2)}{n-4}}|u|^{\frac{2n}{n-4}}dx\Big)^{\frac{n-4}{n}}
\la{ri1-s}
\ee
and
\be
\int_{\Omega}(\Delta u)^2 dx - \frac{n^2}{4}
\int_{\Omega}\frac{|\nabla u|^2}{|x|^2} dx
\geq C_n' \Big(  
\int_{\Omega} X(|x|/D)^{\frac{2(n-1)}{n-2}}|\nabla u|^{\frac{2n}{n-2}}dx\Big)^{\frac{n-2}{n}}
\la{ri2-s}
\ee

The present article contains two main results.
The first theorem extends inequality
\eqref{ri2-s} to dimensions $n=3,4$.
\begin{theorem}
\label{thm2}
Let $\Omega\subset\R^n$, $n=3$ or $n=4$, be a bounded domain and let
$D= \sup_{x \in \Omega} |x|$. There exists $C>0$
such that: \nl
$\ia$ If $n=3$ then
\[
 \int_{\Omega}(\Delta u)^2dx -
\frac{25}{36}
\int_{\Omega}\frac{|\nabla u|^2}{|x|^2}  dx \geq C
\bigg(\int_{\Omega}|\nabla u|^6 X^4(|x|/D) dx\bigg)^{\frac{1}{3}} , 
\quad u \in \cico.
\]
$\ib$ If $n=4$ then
\[
\int_{\Omega}(\Delta u)^2dx -
3\int_{\Omega}\frac{|\nabla u|^2}{|x|^2}  dx \geq C
\bigg(\int_{\Omega}|\nabla u|^4  dx\bigg)^{\frac12}, 
\quad u \in \cico.
\]
Moreover the power $X^4$ in case $n=3$ is the best possible.
\end{theorem}
It is remarkable that in case $n=4$ no logarithmic factor is
required at the RHS, as opposed to the cases $n=3$ and
$n\geq 5$.

Concerning inequality \eqref{ri1-s}, let us first recall what is known
for the corresponding Hardy-Sobolev problem.
In \cite{AFT} it was shown that for any bounded domain
$\Omega\subset\R^n$, $n\geq 3$, and for any $u\in\cico$
there holds
\bean
&& \hspace{-1.5cm} \int_{\Omega}|\nabla u|^2dx -\Big( \frac{n-2}{2}\Big)^2
\int_{\Omega}\frac{u^2}{|x|^2} dx \\
&\geq& (n-2)^{-\frac{2(n-1)}{n}}   S_{1,n}\bigg( \int_{\Omega} X^{\frac{2(n-1)}{n-2}} \big(|x|/D \big) |u|^{\frac{2n}{n-2}}  dx \bigg)^{\frac{n-2}{n}}
\eean
where
\[
S_{1,n} = \pi n(n-2) 
\Big( \frac{\Gamma(\frac{n}{2})}{\Gamma(n)} \Big)^{\frac2n} \ ,
\]
is the best Sobolev constant for the standard Sobolev inequality
in $\R^n$. Moreover the constant $(n-2)^{-\frac{2(n-1)}{n}}S_{1,n}$
is the best possible.
Similarly, in the article \cite{BT19} Sobolev improvements with
best constants were obtained to sharp Hardy
inequalities in Euclidean and hyperbolic space.
We note that by slightly adapting \cite[Theorem 5]{BT19}
we obtain that if $\Omega$ is a bounded domain 
in $\R^n$, $n\geq 3$, then
\bea
&& \int_{\Omega}|\nabla u|^2dx -\Big( \frac{n-2}{2}\Big)^2 
\int_{\Omega}\frac{u^2}{|x|^2}dx
+ \frac{(n-1)(n-3)}{4} \int_{\Omega}\frac{u^2}{|x|^2}
X^2(|x|/D) dx \nonumber \\
& \geq & S_{1,n}
 \Big( \int_{\Omega}   X^{\frac{2(n-1)}{n-2}}
 (|x|/D) 
 |u|^{\frac{2n}{n-2}}  dx \Big)^{\frac{n-2}{n}} 
\la{eli}
\eea
for all $u\in\cic(\Omega)$ and the constant
$S_{1,n}$ is sharp.

The second theorem of this article provides an
estimate with best Sobolev constant for a slightly modified
version of \eqref{ri1-s} which is in the spirit of
\eqref{eli}.
\begin{theorem}
\la{thm1}
Let $\Omega\subset\R^n$, $n\geq 5$, be a bounded domain
and let $D=\sup_{\Omega}|x|$.
For any $u\in\cic (\Omega)$ there holds
\begin{align*}
& \int_{\Omega}(\Delta u)^2 dx - \frac{n^2(n-4)^2}{16}
\int_{\Omega}\frac{u^2}{|x|^4} dx  
 + \frac{n^2(n-4)^2}{16}
 \int_{\Omega}\frac{u^2}{|x|^4} 
X^{\frac{2(n-2)}{n-1}} dx \\[0.1cm]
& \qquad \geq S_{2,n}  \Big(  
\int_{\Omega} X^{\frac{2(n-2)}{n-4}}
|u|^{\frac{2n}{n-4}}dx\Big)^{\frac{n-4}{n}};
\end{align*}
here $X=X(|x|/D)$. Moreover the constant $S_{2,n}$
is the best possible.
\end{theorem}
The proof of Theorem \ref{thm2} is in Section \ref{sec2}
and the proof of Theorem \ref{thm1} is in Section \ref{sec3}.

\section{Rellich-Sobolev inequality I}
\la{sec2}

In this section we shall prove Theorem \ref{thm2}.
An important tool will be the decomposition of functions in 
spherical harmonics \cite[Section IV.2]{StWe}.

We recall that the eigenvalues of the Laplace-Beltrami operator on the unit 
sphere $\sph$ are given by
\[
\mu_k = k(k+n-2) , \qquad k=0,1,2\ldots,
\]
Each $\mu_k$ has multiplicity
\[
d_k ={n+k-1 \choose k} - {n+k-3 \choose k-2}, \quad k\geq 2,
\]
while $d_0=1$ and $d_1=n$.

Let $\{\phi_{kj}\}_{j=1}^{d_k}$ be an orthonormal
basis of eigenfunctions for the eigenvalue $\mu_k$.
Then any function $u\in L^2(\R^n)$ can be decomposed as
\be
u(x)= \sum_{k=0}^{\infty} \sum_{j=1}^{d_k}
u_{kj}(x) =\sum_{k=0}^{\infty} \sum_{j=1}^{d_k}
f_{kj}(r)\phi_{kj}(\omega) 
\la{doublesum}
\ee
where $x=r\omega$, $r>0$, $\omega\in \sph$, and
\[
f_{kj}(r) =\int_{\sph} u(r\omega)\phi_{kj}
(\omega)dS(\omega). 
\]
We note that each $\phi_{kj}$ is the restriction on the 
unit sphere of a harmonic homogeneous polynomial of degree $k$ \cite{StWe}.

Assume now that $u\in\cic(\R^n)$. Since any homogeneous 
polynomial can be written as a linear combination of 
harmonic homogeneous polynomials, taking the Taylor 
expansion of $u$ near the origin we easily infer that
\be
\la{stw}
f_{kj}(r) =O(r^k), \quad
f_{kj}'(r) =O(r^{k-1})
 \; , \qquad \mbox{ as }r\to 0.
\ee
for any $k\geq 1$ and any $j=1,\ldots,d_k$.

We note that
\be
\mu_k \geq n-1 \; , \qquad  \forall \, k\geq 1 , 
\la{mk}
\ee
an estimate that will be used several times
in what follows.

In what follows we shall use
$\sum_{k,j}$ as
a shorthand for  $\sum_{k=0}^{\infty} \sum_{j=1}^{d_k}$.

For simplicity we shall denote by $u_0$
(instead of $u_{01}$) the first (radial) term in the 
decomposition \eqref{doublesum} of $u$ into spherical harmonics.
We note the relation
\be
\int_{\R^n} (\Delta u- \Delta u_0) ^2 dx =
\sum_{k=1}^{\infty}
\sum_{j=1}^{d_k} \int_{\R^n}(\Delta u_{kj})^2dx \, .
\la{la:sph:harm}
\ee

\begin{lemma}
Let $n\geq 3$. For any $u\in C^{\infty}_c(\R^n)$ there holds
\bean
\ia  \;\;  &&  \int_{\R^n}(\Delta u)^2dx = \sum_{k,j} \left\{ \int_0^{\infty}r^{n-1} f_{kj}''^2 dr   \right. \\
&& \left. + (n-1+2\mu_k)  \int_0^{\infty}r^{n-3}
f_{kj}'^{\, 2} dr 
 + \big(2 (n-4) \mu_k +\mu_k^2 \big) \int_0^{\infty}r^{n-5}f_{kj}^2dr
\right\}  \\
\ib  \;\;  &&  \int_{\R^n}
\frac{ |\nabla u|^2}{|x|^2}dx = 
\sum_{k,j}\left\{ \int_0^{\infty}r^{n-3} f_{kj}'^{\, 2} dr
+ \mu_k  \int_0^{\infty}r^{n-5}f_{kj}^2 dr  \right\}  
\eean
\label{lem:g:syn}
\end{lemma}
\proof
Using the orthonormality of the set
$\{\phi_{kj}\}$ we have 
\bean
 \int_{\R^n}(\Delta u)^2dx& =& \sum_{k,j}\int_{\R^n}(\Delta u_{kj})^2dx  \\
& =&  \sum_{k,j}
\int_{0}^{\infty} \Big(  f_{kj}''  + \frac{n-1}{r} f_{kj}'  
 -\frac{\mu_k}{r^2}f_{kj} \Big)^2r^{n-1}dr.
\eean
Part (i) then follows by expanding the square and integrating by parts. Estimates \eqref{stw} ensure that
no terms appear from $r=0$.
The proof of (ii) is similar
and is omitted. 
\endproof

For $n\geq 3$ we set
\[
\mI[u] = \int_{\R^n}(\Delta u)^2 dx - c_n
\int_{\R^n}\frac{|\nabla u|^2}{|x|^2} dx
\]
where the constant $c_n$ is given by \eqref{c_n}.
\begin{lemma}
\la{TZ:thm3.3}
Assume that $n=3$ or $n=4$. There exists $c>0$ such that for 
any $u\in C^{\infty}_c(\R^n)$ there holds
\be
\mI[u] \geq \mI [u_0] + \sum_{j=1}^n \mI [u_{1j}] +
 c \int_{\R^n} \big(\Delta u- \Delta  u_0 -
 \sum_{j=1}^n\Delta  u_{1j} \big)^2 dx .
\la{sy9.52}
\ee
\end{lemma}
\proof
Let $u\in \cic(\R^n)$. Because of the relation
\[
\mI[u] =\mI[u_0]+ \sum_{j=1}^n \mI[u_{1j}] +\sum_{k=2}^{\infty} \sum_{j=1}^{d_k}\mI[u_{kj}],
\]
inequality \eqref{sy9.52} will follow if we establish 
the existence of $c>0$ such that
\be
\mI[u_{kj}] \geq  c \int_{\R^n} (\Delta u_{kj})^2 dx \, , 
\qquad k \geq 2 \; , \;\; 1\leq j\leq d_k.
\la{ploi2}
\ee
Assume first that $n=3$. Let $\lambda>0$ be fixed. For $k\geq 2$ we have
$\mu_k \geq 6$ and therefore
\bean
&& \int_{\R^3} (\Delta u_{kj})^2 dx \\
&=&  \int_0^{\infty}r^2 f_{kj}''^2 dr
+ (2+2\mu_k)  \int_0^{\infty}f_{kj}'^2 dr   + (-2\mu_k +\mu_k^2 ) 
\int_0^{\infty}\!\! r^{-2}f_{kj}^2dr    \\
&\geq&  \Big( \frac{9}{4} +2\lambda \mu_k  \Big) 
\int_0^{\infty}\!\! 
f_{kj}'^2 dr  + \Big(  2(1-\lambda)\frac{1}{4}\mu_k -2\mu_k +\mu_k^2  \Big)  \int_0^{\infty}\!\! r^{-2}f_{kj}^2dr    \\
&\geq&  \Big(  \frac{9}{4} +12\lambda \Big) \int_0^{\infty}
f_{kj}'^2 dr  + \Big(  \frac{9}{2} -\frac{\lambda}{2} \Big)
\mu_k  \int_0^{\infty}r^{-2}f_{kj}^2dr  .
\eean
Choosing $\lambda =9/50$ we arrive at
\[
\int_{\R^3} (\Delta u_{kj})^2 dx \geq \frac{441}{100}
\int_{\R^3}\frac{|\nabla u_{kj}|^2}{|x|^2}dx ,
\]
and \eqref{ploi2} follows. 
In case $n=4$ we argue similarly. We now have
$\mu_k\geq 8$, hence
\bean
 \int_{\R^4} (\Delta u_{kj})^2 dx 
&=&  \int_0^{\infty}r^3 f_{kj}''^2 dr
+ (3+2\mu_k)  \int_0^{\infty}r f_{kj}'^2 dr +\mu_k^2 
\int_0^{\infty}\!\! r^{-1}f_{kj}^2dr    \\
&\geq&  (4+2\mu_k)  \int_0^{\infty}r f_{kj}'^2 dr
+\mu_k^2 
\int_0^{\infty}\!\! r^{-1}f_{kj}^2dr   \\
& \geq  & 8\int_{\R^4}\frac{|\nabla u_{kj}|^2}{|x|^2}dx,
\eean
as required. 
\endproof

\begin{lemma}
\la{TZ:lem3.4}
Let $n=3$ or $n=4$. Then there exists $c>0$ such that
\be
\la{sy10.1}
\mI[u_0] \geq c \bigg( \int_{B_1} |\nabla u_0|^{\frac{2n}{n-2}}  dx  \bigg)^{\frac{n-2}{n}}.
\ee
Additionally for $n=3$ we have
\be
\la{sy10.2}
\mI[u_{1j}] \geq c \bigg( \int_{B_1} |\nabla u_{1j}|^6 X^4 dx  \bigg)^{\frac{1}{3}} \!\! , \qquad j=1,2,3,
\ee
while for $n=4$
\be
\la{sy10.3}
\mI[u_{1j}] \geq c \bigg( \int_{B_1} |\nabla u_{1j}|^4  dx  
\bigg)^{\frac{1}{2}}\!\!  , \qquad  j=1,2,3,4.
\ee
Here $X=X(|x|)$.
\end{lemma}
\proof
From  Lemma \ref{lem:g:syn} (i) and the standard Sobolev inequality we obtain
\bean
\mI[u_0]  & \geq &  \int_0^1 f_0''^2 \, r^{n-1}dr   \\
&\geq& c \bigg( \int_0^1 |f_0'|^{\frac{2n}{n-2}} r^{n-1}dr \bigg)^{\frac{n-2}{n}}\\
&=& c \bigg( \int_{B_1} |\nabla u_0|^{\frac{2n}{n-2}}  dx  \bigg)^{\frac{n-2}{n}}
\eean
as required. 

Assume now that $n=3$. By Lemma \ref{lem:g:syn}
and the improved Hardy-Sobolev inequality of
\cite{AFT} we have
\bean
\mI[u_{1j}]  & =&  \int_0^1 f_{1j}''^2 r^2 dr 
- \frac{1}{4}  \int_0^1 f_{1j}'^2 dr  \\
&& \hspace{1.5cm}  + \frac{50}{9} \Big( \int_0^1 f_{1j}'^2 dr -
\frac14 \int_0^1 r^{-2}f_{1j}^2 dr \Big)  \\
&\geq& c \bigg( \int_0^1 |f_{1j}'|^6
X^4 \, r^2 dr \bigg)^{\frac{1}3} 
+ c \bigg( \int_0^1 |f_{1j}|^6
X^4 \, dr \bigg)^{\frac{1}3} \\
&\geq& c \bigg( \int_{B_1} |\nabla u_{1j}|^6
X^4 dx  \bigg)^{\frac{1}{3}} .
\eean
In case $n=4$ we argue similarly applying 
again Lemma \ref{lem:g:syn} and, now,
the standard Sobolev inequality; we obtain
\bean
\mI[u_{1j}] &=& 
\int_0^1 f_{1j}''^2 \, r^3 dr + 6 \int_0^1 f_{1j}'^2 \, r\,  dr \\
&\geq & c \bigg( \int_0^1 |f_{1j}'|^4   \, r^3 dr \bigg)^{\frac12} 
+  c \bigg( \int_0^1 |f_{1j}|^4   \, 
r \, dr \bigg)^{\frac12} \\ 
&\geq & c \bigg( \int_{B_1} |\nabla u_{1j}|^4  dx \bigg)^{\frac12} ,
\eean
as required. 
\endproof

\noindent
{\bf \em Proof of Theorem \ref{thm2}.} We first note that by the standard Sobolev inequality
we have
\[
\ino  (\Delta u-  \Delta u_0 - 
\sum_{j=1}^n\Delta u_{1j})^2 dx \geq
c \bigg( \ino  | \nabla u- \nabla u_0 - \sum_{j=1}^n
\nabla u_{1j}|^{\frac{2n}{n-2}}
\, dx \bigg)^{\frac{1}{3}} \!\! ;
\]
In case $n=3$ we apply \eqref{sy9.52},
\eqref{sy10.1}, \eqref{sy10.2} and the triangle
inequality to obtain
\bean
\mI[u] 
 &\geq& \mI [u_0] + \sum_{j=1}^n \mI [u_{1j}] +
 c \int_{\R^n} \big(\Delta u- \Delta  u_0 -
 \sum_{j=1}^n\Delta  u_{1j} \big)^2 dx \\
 &\geq& c \bigg ( \ino  | \nabla u_0|^6 X^4 dx  \bigg)^{\frac{1}{3}}  
 +c\sum_{j=1}^n  \bigg( \int_{B_1} |\nabla u_{1j}|^6 X^4 dx  \bigg)^{\frac{1}{3}} \\
 && + c \bigg( \ino  |\nabla u-\nabla u_0
- \sum_{j=1}^n\nabla u_{1j} |^6  dx 
\bigg)^{\frac{1}{3}} \\
&\geq& c \bigg ( \ino |\nabla u|^6
X^4 dx  \bigg)^{\frac{1}{3}}.
\eean
In case $n=4$ we argue similarly, the only
difference being that we use \eqref{sy10.3}
instead of \eqref{sy10.2}.

We next prove the optimality of the power $X^4$ in
(i), that is in case $n=3$. So let us assume instead that there exist $\mu<4$ and $c>0$ so that
\be
 \int_{\Omega}(\Delta u)^2dx -
\frac{25}{36}
\int_{\Omega}\frac{|\nabla u|^2}{|x|^2}  dx \geq c
\bigg(\int_{\Omega}|\nabla u|^6 X^{\mu}(|x|/D) dx\bigg)^{\frac{1}{3}} , 
\la{end}
\ee
for all $u \in \cico$.
Without loss of generality we assume that
$B_1\subset\Omega$.
We consider small positive
numbers $\epsilon$ and $\delta$
and define the functions
\[
u_{\epsilon,\delta}(x) =
f_{\epsilon,\delta}(r)\phi_1(\omega) :=
r^{\frac{1}{2} +\epsilon}
X(r)^{-\frac{1}{2} +\delta} \psi(r) \phi_1(\omega)
\]
where $\phi_1(\omega)$ is a normalized eigenfunction 
for the first non-zero eigenvalue of the
Laplace-Beltrami operator on ${\rm S}^2$ 
and $\psi(r)$ is a smooth radially symmetric function 
supported in $B_1$ and equal to one near $r=0$.

Applying Lemma \ref{lem:g:syn} we see that
$\int(\Delta u_{\epsilon,\delta})^2dx -\frac{25}{36}
\int \frac{ |\nabla u_{\epsilon,\delta}|^2}{|x|^2}dx$
is a linear combination of the integrals
\[
I_{\epsilon,\delta}^{(j)}=
\int_0^1r^{-1+2\epsilon}X^{-1+j+2\delta}\psi^2 dr, \quad
0\leq j\leq 4,
\]
and of integrals
that contain at least one derivative
of $\psi$ and are, therefore, uniformly bounded.
Moreover simple computations yield that
for $j=3,4$ the integrals
$I_{\epsilon,\delta}^{(j)}$ 
are also uniformly bounded for small $\epsilon,\delta>0$.

Restricting attention to a small neighbourhood of the origin
where $\psi=1$ we find
\[
f_{\epsilon,\delta}'(r) =r^{-\frac12 + \epsilon}\Big(
\big(\frac12 +\epsilon\big)X^{-\frac12 +\delta} +
\big(-\frac12 +\delta\big)
X^{\frac12 +\delta}  \Big)
\]
and
\[
f_{\epsilon,\delta}''(r) =r^{-\frac32 +\epsilon}
\Big(
\big(\epsilon^2 -\frac14 \big) X^{-\frac12 +\delta} +2\epsilon
\big(-\frac12 +\delta\big) X^{\frac12 +\delta}
+  \big(\delta^2 -\frac14\big)X^{\frac32 +\delta}
\Big)
\]
Hence we arrive at
\bean
&& \hspace{-1.2cm}\int_{B_1}(\Delta u_{\epsilon,\delta})^2dx -\frac{25}{36}
\int_{B_1} \frac{ |\nabla u_{\epsilon,\delta}|^2}{|x|^2}dx \\
&=&\Big(  \frac{191}{36}\epsilon +
\frac{173}{36}\epsilon^2 +\epsilon^4\Big)
I_{\epsilon,\delta}^{(0)} \\
&& - \Big( \frac{191}{72} - \frac{191}{36}\delta 
+\big( \frac{173}{36} - \frac{173}{18} \delta\big)
\epsilon  +(2-4\delta )\epsilon^3\Big)
I_{\epsilon,\delta}^{(1)} \\
&&+ \Big( \frac{209}{144} -\frac{191}{36}\delta
+\frac{173}{36}\delta^2 + \big(  
\frac12  -4\delta +6\delta^2 \big)\epsilon^2\Big)
I_{\epsilon,\delta}^{(2)}
+O(1).
\eean
It is easily seen that
\[
I_{\epsilon,0}^{(j)} =\frac{1}{2\epsilon} +O(1) \; , \quad  j=0,1,2.
\]
Hence, rearranging also terms we obtain
\bean
\int_{B_1}(\Delta u_{\epsilon,\delta})^2dx -\frac{25}{36}
\int_{B_1} \frac{ |\nabla u_{\epsilon,\delta}|^2}{|x|^2}dx 
&=&  \frac{191}{72} \big( 2\epsilon I_{\epsilon,\delta}^{(0)}
-(1-2\delta) I_{\epsilon,\delta}^{(1)} \big) \\
&& + \Big( \frac{209}{144} -\frac{191}{36}\delta
+\frac{173}{36}\delta^2 \Big) I_{\epsilon,\delta}^{(2)} +O(1).
\eean
Now, by \cite[p181]{BFT} we have
\[
2\epsilon I_{\epsilon,\delta}^{(0)} -(1-2\delta)
I_{\epsilon,\delta}^{(1)} =O(1).
\]
Hence, letting $\epsilon\to 0$ we obtain
\bean
 \int_{B_1}(\Delta u_{\epsilon,\delta})^2dx -\frac{25}{36}
\int_{B_1} \frac{ |\nabla u_{\epsilon,\delta}|^2}{|x|^2}dx  
&\to & \Big( \frac{209}{144} -\frac{191}{36}\delta
+\frac{173}{36}\delta^2 \Big)
I_{0,\delta}^{(2)} +O(1) \\
&=& \frac{209}{144} \int_0^1r^{-1}X^{1+2\delta}\psi^2 dr   +O(1),
\eean
which is finite for any $\delta>0$ and diverges to infinity as
$\delta\to 0+$.

Now, for $\delta> (4-\mu)/6$ we have
\[
\int_{B_1}|\nabla u_{\epsilon,\delta}|^6 X^{\mu} dx
\geq c\int_0^{1/2} r^{-1+6\epsilon}X^{\mu-3+6\delta}dr.
\]
Letting first $\epsilon\to 0$ and then $\delta \to
\frac{4-\mu}{6}+ $ the last integral tends to infinity.
Hence the Rayleigh quotient tends to zero, which
implies that the constant $c$ in \eqref{end} should be zero. This concludes the proof.
$\hfill\Box$

\section{Rellich-Sobolev inequality II}
\la{sec3}

In this section we shall prove Theorem \ref{thm1}. Throughout the proof
we shall make use of spherical coordinates $(r,\omega)$, $r>0$, $\omega\in\sph$.
We denote 
by $\nabla_{\omega}$ and $\Delta_{\omega}$ the gradient and Laplacian on
$\sph$.
\begin{lemma}
\la{lem:sph_lapl}
Let $\theta\in\R$. For any $v\in\cic(\R^n\setminus\{0\})$ there holds
\bean
&& \hspace{-.7cm}\int_{\R^n}(\Delta v)^2 |x|^{\theta}dx \\
&=& \int_0^{\infty} \!\!\int_{\sph} v_{rr}^2 r^{n+\theta-1} dS \, dr +
(n-1)(1-\theta)\int_0^{\infty}\!\! \int_{\sph} v_r^2 r^{n+\theta-3} dS \, dr
\nonumber \\  
&& + \int_0^{\infty} \!\!
\int_{\sph} (\Delta_{\omega}v)^2 r^{n+\theta -5} dS \, dr
+2 \int_0^{\infty} \!\! \int_{\sph} |\nabla_{\omega}v_r|^2 r^{n+\theta -3} \, dS \, dr 
\nonumber \\
&& -(\theta-2)(n+\theta -4)  \int_0^{\infty}\!\! \int_{\sph} 
|\nabla_{\omega}v|^2 r^{n+\theta -5} dS \, dr \, .
\eean
\end{lemma}
\proof This follows by writing
\[
\Delta v = v_{rr} +\frac{n-1}{r}v_r +\frac{1}{r^2}\Delta_{\omega}v 
\]
and integrating by parts; we omit the details.
\endproof

In the next lemma and also later, we shall use subscripts R and NR to denote
the radial and non-radial component of a given functional.
\begin{lemma}
\la{lem1}
Let $n\geq 5$, $\beta>0$ and define
\[
A= \frac{1}{\beta^2}\big( 2n-4-\beta(n-4+\beta)\big).
\]
Let $u\in\cic(\R^n)$.
Changing variables by  $u(r,\omega) =y(t,\omega)$, $t=r^{\beta}$, we have
\[
\frac{\Int_{\R^n}(\Delta u)^2dx}{\Big( \Int_{\R^n}|u|^{\frac{2n}{n-4}}dx 
\Big)^{\frac{n-4}{n}}}
= \beta^{\frac{4(n-1)}{n}} \frac{  \cA_{\rm R}[y]  + \cA_{\rm NR}[y]}
{ \Big( \Int_0^{\infty}\!\!\Int_{\sph} t^{\frac{n-\beta}{\beta}} 
 |y|^{\frac{2n}{n-4}}dS \, dt  \Big)^{  \frac{n-4}{n}}}
\]
where
\bean
\cA_{\rm R}[y]& =& \int_0^{\infty}\!\! \int_{\sph}
 \big( t^{\frac{3\beta+n-4}{\beta}} y_{tt}^2
+A t^{\frac{\beta+n-4}{\beta}} y_{t}^2 \big) dS \, dt \\
\cA_{\rm NR}[y]& =& \int_0^{\infty}\!\! \int_{\sph}
\Big(\frac{1}{\beta^4} t^{\frac{n-\beta-4}{\beta}}
 (\Delta_{\omega}y)^2  + \frac{2}{\beta^2}
t^{\frac{\beta+n-4}{\beta}}|\nabla_{\omega}y_t|^2 \\
&& \hspace{2cm} +\frac{2(n-4)}{\beta^4} 
t^{\frac{n-\beta-4}{\beta}}
|\nabla_{\omega}y|^2 \Big) dS \, dt
\eean
\end{lemma}
\proof 
After some lengthy but otherwise elementary computations we find
\[
\int_0^{\infty} \big( u_{rr} + \frac{n-1}{r}u_r \big)^2 r^{n-1} dr = \beta^3 \int_0^{\infty} \big( t^{\frac{3\beta+n-4}{\beta}} y_{tt}^2
+A t^{\frac{\beta+n-4}{\beta}} y_{t}^2 \big) dt
\]
and
\[
\int_0^{\infty} |u|^{\frac{2n}{n-4}} r^{n-1} dr
= \frac{1}{\beta}\int_0^{\infty} |y|^{\frac{2n}{n-4}} t^{\frac{n-\beta}{\beta}} 
 dt \, .
\]
Similar computations involving the non-radial (tangential) 
derivatives yield the term $\cA_{\rm NR}[y]$. We omit the
details. \endproof

We now consider the Rayleigh quotient  for the Rellich-Sobolev inequality
\eqref{ri1-s}.
Changing variables by $u(x) =|x|^{-\frac{n-4}{2}}v(x)$ we obtain (cf. \cite[Lemma 2.3 (ii)]{TZ})
\bea
&& \int_{\Omega} (\Delta u)^2 dx
 -\frac{ n^2(n-4)^2}{16} \int_{\Omega} \frac{u^2}{|x|^4}dx  
 \la{ippp} \\
&=& \Int_{\Omega} 
\!\!\Big( |x|^{4-n}(\Delta v)^2   + \frac{n(n-4)}{2}  |x|^{2-n}|\nabla v|^2 
 -n(n-4)  |x|^{-n}(x\cdot \nabla v)^2 \Big) dx \, . \nonumber \\
&=:& J[v]
\eea
Applying Lemma \ref{lem:sph_lapl} we find that
\bea
J[v]&=& \int_0^1 \!\!\int_{\sph} r^3v_{rr}^2  dS \, dr +
\frac{n^2-4n+6}{2} \int_0^1 \!\!\int_{\sph} rv_r^2 dS \, dr
\nonumber \\  
&& + \int_0^1 \!\!\int_{\sph} r^{-1}(\Delta_{\omega}v)^2 dS \, dr
+2 \int_0^1 \!\!\int_{\sph} |\nabla_{\omega}v_r|^2 r \, dS \, dr 
\nonumber \\
&& + \frac{n(n-4)}{2} \int_0^1 \!\! \int_{\sph} 
r^{-1}|\nabla_{\omega}v|^2 dS \, dr \, .
\la{ipp5}
\eea
In view of \eqref{ipp5} we set
\bea
J_{\rm R}[v]&= & \int_0^1 \!\!\int_{\sph}   r^3v_{rr}^2  dS \, dr +
\frac{n^2-4n+6}{2} \int_0^1 \!\!\int_{\sph}  rv_r^2 dS \, dr \nonumber \\
J_{\rm NR}[v] &=&  \int_0^1 \!\!\int_{\sph} r^{-1}(\Delta_{\omega}v)^2 dS \, dr
+2 \int_0^1 \!\!\int_{\sph} r \, |\nabla_{\omega}v_r|^2  dS \, dr 
\nonumber \\
&& + \frac{n(n-4)}{2} \int_0^1 \!\!\int_{\sph} 
r^{-1}|\nabla_{\omega}v|^2 dS \, dr \nonumber  ,
\la{3}
\eea
the radial and non-radial parts of $J[v]$, so that,
\[
J[v] =  J_{\rm R}[v]+ J_{\rm NR} [v].
\]
We shall change variables once more and for this
we define the functions
\be
g(r)= \exp\Big( 1-X(r)^{- \frac{n}{2(n-1)}}\Big) \; , \qquad
\alpha(r) =X(r)^{-\frac{3(n-2)}{4(n-1)}}
g(r)^{ \frac{n-4}{2\beta}}.
\la{gr}
\ee
\begin{lemma}
Let $n\geq 5$, $\beta>0$ and set
\[
s=\frac{n-4}{2\beta}.
\]
Let $v\in\cic (B_1\setminus\{0\})$. Changing variables by
\be
v(r,\omega) =\alpha(r) w(t,\omega)  \; , \qquad t=g(r) ,
\la{5}
\ee
we have
\begin{align*}
\ia \quad & J_{\rm R}[v] =  \int_0^1\!\!\int_{\sph} \Big\{   \Big( \frac{n}{2(n-1)}
\Big)^3  t^{\frac{3\beta +n-4}{\beta}}
w_{tt}^2  + t^{\frac{\beta +n-4}{\beta}}G(t) w_{t}^2   \\
& \hspace{5cm} +  t^{\frac{-\beta +n-4}{\beta}}H(t) w^2 \Big\} dS \, dt  \\
\ib \quad &  J_{\rm NR}[v] = 
\frac{2(n-1)}{n} \int_0^1 \!\!\int_{\sph} t^{ \frac{n-\beta-4}{\beta}}
 X(t)^{\frac{8-4n}{n}}  (\Delta_{\omega}w)^2 dS \, dt \\
 & \quad\qquad\qquad  + \frac{n}{n-1}  \int_0^1 \!\!\int_{\sph} t^{ \frac{n+\beta-4}{\beta}}
 X(t)^{\frac{4-2n}{n}}  |\nabla_{\omega}w_t|^2 dS \, dt  \\
 & \quad\qquad\qquad + \int_0^1 \!\!\int_{\sph} t^{ \frac{n-\beta-4}{\beta}} |\nabla_{\omega}w|^2 
 K(t) dS \, dt \\
\ic \quad &  \int_0^1 \!\!\int_{\sph}
r^{-1}X(r)^{ \frac{2n-4}{n-4}}|v|^{\frac{2n}{n-4}}dS \, dr
= \frac{2(n-1)}{n}\int_0^1 \!\!\int_{\sph} |w|^{\frac{2n}{n-4}} t^{ \frac{n-\beta}{\beta}} dS \, dt \, ,
\end{align*}
where the functions $G(t)$, $H(t)$ and $K(t)$ are given by
\begin{align*}
G(t) &= \frac{ n(n^2-4n+8)}{4(n-1)} X(t)^{\frac{4-2n}{n}} \!
 -\frac{n^3(2s^2+2s+1)}{8(n-1)^3} 
 + \frac{5n(n-2)(3n-2)}{16(n-1)^3} X(t)^2  \\
H(t) &=  -\frac{s^2n(n^2-4n+8)}{4(n-1)}
X(t)^{\frac{4-2n}{n}}
 +\frac{s(n-2)(n^2-4n+8)}{2(n-1)} X(t)^{\frac{4-n}{n}} \\
& \qquad +\frac{s^4n^3}{8(n-1)^3} 
+\frac{3(n^2-4)(n^2-4n+8)}{16n(n-1)}X(t)^{\frac{4}{n}}  \\
& \qquad - \frac{5s^2n(n-2)(3n-2)}{16(n-1)^3}X(t)^2
 - \frac{5sn(n-2)(3n-2)}{8(n-1)^3}X(t)^3  \\
& \qquad -\frac{9(3n-2)(5n-2)(n^2-4)}{128n(n-1)^3}X(t)^4 \\
K(t) &=  (n-1)(n-4)X(t)^{\frac{8-4n}{n}} -\frac{n(n-4)^2}{4(n-1)\beta^2}
X(t)^{\frac{4-2n}{n}} \\
 & \qquad + \frac{(n-2)(n-4)}{(n-1)\beta} X(t)^{\frac{4-n}{n}}
+\frac{3(n^2-4)}{4n(n-1)}X(t)^{\frac{4}{n}}.
\end{align*}
\la{lem:mark}
\end{lemma}
\proof
To prove (i) we set for simplicity
\[
J^*_{\rm R}[v]=  \int_0^1    r^3v_{rr}^2   dr +
\frac{n^2-4n+6}{2} \int_0^1   rv_r^2  dr \, .
\]
We first note that $r$ and $t=g(r)$ are also related by the relation
\be
X(t) =X(r)^{\frac{n}{2(n-1)}}
\la{4}
\ee
and that
\[
dt = \frac{n}{2(n-1)} \frac{g(r)}{r} 
X(r)^{\frac{n-2}{2(n-1)}}dr \, .
\]
Expressing $J^*_{\rm R}[v]$ in terms of the function
$w(t)$ involves some lengthy computations, of which we include
only the main steps.

From \eqref{5} we have
\begin{align*}
v_r =& \; \alpha g' w_t +\alpha' w  \\
v_{rr} =&  \; \alpha g'^2 w_{tt} + (2\alpha' g' +\alpha g'')w_t
+\alpha'' w \, .
\end{align*}
Substuting in $J^*_{\rm R}[v]$ and expanding we find that
\bea
J^*_{\rm R}[v]&=& \big( \frac{n}{2(n-1)} \big)^3 \int_0^1  t^{\frac{3\beta +n-4}{\beta}}
w_{tt}^2 dt +\int_0^1 B(t)w_{t}^2 dt 
+\int_0^1 C(t)w^2 dt \nonumber \\
&& + \int_0^1 D(t)w_{tt} \, w_{t} dt  +
 \int_0^1 E(t)w_{tt}w \, dt  + \int_0^1 F(t)w_t \, w \, dt 
\la{epan}
\eea
where the functions $B(t),\ldots,F(t)$ will be described
below in terms of the variable $r$. 
Integrating by parts we obtain from \eqref{epan} that
\[
J^*_{\rm R}[v]= \big( \frac{n}{2(n-1)} \big)^3 \int_0^1  t^{\frac{3\beta +n-4}{\beta}}
w_{tt}^2 dt +\int_0^1 P(t)w_{t}^2 dt 
+\int_0^1 Q(t)w^2 dt 
\]
where
\be
\la{7}
P(t) =B(t) -\frac{1}{2}D_t(t) -E(t)  \; , \qquad Q(t) =C(t) +\frac12 E_{tt}(t) -
\frac12 F_t(t) \, .
\ee
To compute the functions $P(t)$ and $Q(t)$
it is convenient to regard them
as functions of the variable $r$. To do this we consider
the functions $B,C,D,E$  and $F$ also as functions of $r$ and
indicate this with tildes; we shall thus write
$B(t)=\tilde{B}(r)$, etc.
Relations \eqref{7} then take the form
\be
\tilde{P}(r) =\tilde{B} -\frac{1}{2g'}\tilde{D}_r -\tilde{E} \; , \qquad
\tilde{Q}(r)=\tilde{C}+ \frac12 \Big(\frac{\tilde{E}_{rr}}{g'^2} -\frac{ g''\tilde{E}_r}{g'^3}\Big) -\frac{1}{2g'}\tilde{F}_r \, .
\la{mc}
\ee
After some computations we eventually find
\begin{align*}
\tilde{B}(r)&= \frac{r^3}{g'}\big( 2\alpha'g' +\frac{n-1}{r}\alpha g'  +\alpha g'' \big)^2 -\frac{n(n-4)}{2} r\alpha^2 g' \\
\tilde{C}(r) &= \frac{r^3}{g'}\big( \alpha'' +\frac{n-1}{r}\alpha'\big)^2 -\frac{n(n-4)}{2} \frac{r}{g'}\alpha'^2 \\
\tilde{D}(r) &= 2r^3\alpha g' \big(  2\alpha'g' +\frac{n-1}{r}\alpha g'  +\alpha g''  \big) \\
\tilde{E}(r) &=  2r^3\alpha g' \big( \alpha'' +\frac{n-1}{r}\alpha'\big) \\
\tilde{F}(r) &=  2r^3 \big( 2\alpha' + \frac{n-1}{r}\alpha + 
\alpha\frac{g''}{g'} \big)\big( \alpha'' +\frac{n-1}{r}\alpha'\big) - n(n-4)r\alpha \alpha' \, .
\end{align*}
Substituting in \eqref{mc} we arrive at
\begin{align*}
\tilde{P}(r)&= 2r^3\alpha'^2 g' -6r^2\alpha\alpha' g'  + \frac{n^2-4n+6}{2}
r\alpha^2 g'   -3r^2\alpha^2g'' \\[0.2cm]
& \qquad  -4r^3\alpha \alpha'' g' -2r^3\alpha\alpha' g'' -r^3\alpha^2 g'''  \\
\tilde{Q}(r)&= \frac{1}{g'}\Big( 6r^2\alpha\alpha''' - \frac{n^2-4n-6}{2}r\alpha\alpha''
 - \frac{n^2-4n+6}{2} \alpha\alpha' +r^3 \alpha\alpha^{(4)}\Big).
\end{align*}
Now, some more computations give
\begin{align*}
g'(r)& =\frac{n}{2(n-1)}  \frac{g(r)}{r} X(r)^{\frac{n-2}{2(n-1)}} \; , \\
g''(r) & = \Big(  -\frac{n}{2(n-1)} X^{\frac{n-2}{2(n-1)}} 
+\frac{n^2}{4(n-1)^2} X(r)^{\frac{n-2}{n-1}} + \frac{n(n-2)}{4(n-1)^2}
X(r)^{\frac{3n-4}{2(n-1)}}\Big)
 \frac{g(r)}{r^2} \\
g'''(r) & = \Big(  -\frac{3n(n-2)}{4(n-1)^2} X(r)^{\frac{3n-4}{2(n-1)}}
+\frac{3n^2(n-2)}{8(n-1)^3} X^{\frac{2n-3}{n-1}}
+
\frac{n(n-2)(3n-4)}{8(n-1)^3}X^{\frac{5n-6}{2(n-1)}}  \\
& \qquad 
+\frac{n}{n-1} X(r)^{\frac{n-2}{2(n-1)}} 
-\frac{3n^2}{4(n-1)^2} X(r)^{\frac{n-2}{n-1}} 
+\frac{n^3}{8(n-1)^3} X(r)^{\frac{3n-6}{2(n-1)} }\Big)
 \frac{g(r)}{r^3}.
\end{align*}
Moreover, 
\begin{align}
\alpha'(r) & = \frac{ g(r)^{ s }}{r}
\Big( \frac{s}{2(n-1)}X^{\frac{2-n}{4(n-1)}} -
\frac{3(n-2)}{4(n-1)} X(r)^{\frac{n+2}{4(n-1)}}  \Big) \nonumber \\
\alpha''(r) &= \frac{ g(r)^{s}  }{r^2} 
\Big( -\frac{sn}{2(n-1)} X(r)^{\frac{2-n}{4(n-1)}} 
 + \frac{s^2n^2}{4(n-1)^2}X(r)^{\frac{n-2}{4(n-1)}}  \nonumber \\
& \quad  + \frac{3(n-2)}{4(n-1)} X(r)^{\frac{n+2}{4(n-1)}}
  \! 
  - \!\frac{sn(n-2)}{2(n-1)^2} X(r)^{\frac{3n-2}{4(n-1)}}
-\!\frac{3(n^2-4)}{16(n-1)^2} X(r)^{\frac{5n-2}{4(n-1)}}\!
\Big) \nonumber \\
\alpha'''(r) &= \frac{ g(r)^{s}  }{r^3} \Big(
 \frac{sn}{n-1}X^{\frac{2-n}{4(n-1)}}
 -\frac{3s^2n^2}{4(n-1)^2}X(r)^{\frac{n-2}{4(n-1)}} \nonumber \\
& \quad - \frac{3(n-2)}{2(n-1)}X(r)^{\frac{n+2}{4(n-1)}}
 + \frac{s^3n^3}{8(n-1)^3}X^{\frac{3n-6}{4(n-1)}} 
 +\frac{3sn(n-2)}{2(n-1)^2}X^{\frac{3n-2}{4(n-1)}} \nonumber \\
&\quad -\frac{3s^2n^2(n-2)}{16(n-1)^3}X(r)^{\frac{5n-6}{4(n-1)}}
  +\frac{9(n^2-4)}{16(n-1)^2}X(r)^{\frac{5n-2}{4(n-1)}} \nonumber \\
&\quad -\frac{sn(n-2)(15n-2)}{32(n-1)^3}X(r)^{\frac{7n-6}{4(n-1)}}
  -\frac{3(n^2-4)(5n-2)}{64(n-1)^3}X(r)^{\frac{9n-6}{4(n-1)}}
\Big)
\la{alpha}
\end{align}
and
\begin{align*}
 \alpha^{(4)} (r) 
& = \frac{ g(r)^{s}}{r^4}
\Big( \frac{3sn}{n-1}X(r)^{\frac{2-n}{4(n-1)}} 
  -  \frac{11s^2n^2}{4(n-1)^2}X(r)^{\frac{n-2}{4(n-1)}}
\\
&\qquad  -\frac{9(n-2)}{2(n-1)}X(r)^{\frac{n+2}{4(n-1)}} 
+ \frac{3s^3n^3}{4(n-1)^3}X(r)^{\frac{3n-6}{4(n-1)}}
  \\
& \qquad +\frac{11sn(n-2)}{2(n-1)^2}X(r)^{\frac{3n-2}{4(n-1)}}
-\frac{ s^4n^4}{16(n-1)^4} X(r)^{\frac{5n-10}{4(n-1)}} \\
 & \qquad  - \frac{ 9s^2n^2(n-2)}{8(n-1)^3}X(r)^{\frac{5n-6}{4(n-1)}} 
 +\frac{33(n^2-4)}{16(n-1)^2}X(r)^{\frac{5n-2}{4(n-1)}} \\
 & \qquad  -\frac{3sn(n-2)(15n-2)}{16(n-1)^3}X(r)^{\frac{7n-6}{4(n-1)}}
   +\frac{5s^2n^2(n-2)(3n-2)}{32(n-1)^4}X(r)^{\frac{9n-10}{4(n-1)}} \\
 & \qquad    - \frac{ 9(5n-2)(n^2-4)}{32(n-1)^3}
 X(r)^{\frac{9n-6}{4(n-1)}}
  +\frac{ 5sn^2(n-2)(3n-2)}{16(n-1)^4}X(r)^{\frac{11n-10}{4(n-1)}} \\
 & \qquad  +\frac{ 9(3n-2)(5n-2)(n^2-4)}{256(n-1)^4}X(r)^{\frac{13n-10}{4(n-1)}}
 \Big).
\end{align*}
Combining the above we eventually arrive at
\bean
\tilde{P}(r)&=& g(r)^{\frac{\beta+n-4}{\beta}} \Big(  
\frac{ n(n^2-4n+8)}{4(n-1)} X(r)^{\frac{2-n}{n-1}}
 -\frac{n^3(2s^2+2s+1)}{8(n-1)^3} \\
 && \hspace{2.5cm} + \frac{5n(n-2)(3n-2)}{16(n-1)^3} X(r)^{\frac{n}{n-1}} \Big)
\eean
and
\begin{align*}
\tilde{Q}(r) & =  g(r)^{\frac{-\beta+n-4}{\beta}}
 \Big( -\frac{s^2n(n^2-4n+8)}{4(n-1)}
X(r)^{\frac{2-n}{n-1}} 
+\frac{s(n-2)(n^2-4n+8)}{2(n-1)} X(r)^{\frac{4-n}{2(n-1)}} \\
& \qquad +\frac{s^4n^3}{8(n-1)^3} 
+\frac{3(n^2-4)(n^2-4n+8)}{16n(n-1)}X(r)^{\frac{2}{n-1}}  \\
& \qquad - \frac{5s^2n(n-2)(3n-2)}{16(n-1)^3}X(r)^{\frac{n}{n-1}} 
- \frac{5sn(n-2)(3n-2)}{8(n-1)^3}X(r)^{\frac{3n}{2(n-1)}}  \\
& \qquad -\frac{9(3n-2)(5n-2)(n^2-4)}{128n(n-1)^3}
X(r)^{\frac{2n}{n-1}}\Big).
\end{align*}
Part (i) now follows by recalling \eqref{4} and noting that
\[
P(t) = t^{ \frac{\beta + n-4}{\beta}}G(t) \; , \qquad
Q(t) = t^{ \frac{-\beta + n-4}{\beta}}H(t) .
\]
To prove part (ii) we first note that
\bean
 \int_0^1 \!\!\int_{\sph} r^{-1}(\Delta_{\omega}v)^2 dS \, dr
 &=&   \int_0^1 \!\!\int_{\sph} r^{-1}\alpha(r)^2 (\Delta_{\omega}w)^2 
 \frac{1}{g'(r)}dS \, dt \\
 &=& \frac{2(n-1)}{n} \int_0^1 \!\!\int_{\sph} t^{ \frac{n-\beta-4}{\beta}}
 X(t)^{\frac{8-4n}{n}}  (\Delta_{\omega}w)^2 dS \, dt 
\eean
and similarly
\[
 \int_0^1 \!\!\int_{\sph} r^{-1}|\nabla_{\omega}v|^2 dS \, dr =   
 \frac{2(n-1)}{n} \int_0^1 \!\!\int_{\sph} t^{ \frac{n-\beta-4}{\beta}}
 X(t)^{\frac{8-4n}{n}}  |\nabla_{\omega}w|^2 dS \, dt  \, .
\]
For the remaining term in $J_{\rm NR}[v]$ we compute
\bean
&& \int_0^1 \!\! \int_{\sph} r |\nabla_{\omega}v_r|^2  dS \, dr \\
&=& \int_0^1 \!\!\int_{\sph} r \alpha^2 g'  |\nabla_{\omega}w_t|^2 \, dS \, dt
-  \int_0^1 \!\! \int_{\sph} |\nabla_{\omega}w|^2 
\frac{1}{g'}(\alpha\alpha'' r+\alpha\alpha')\, dS \, dt
\eean
On the one hand we have
\[
\int_0^1 \!\! \int_{\sph} \alpha^2 g' r |\nabla_{\omega}w_t|^2 \, dS \, dt
=\frac{n}{2(n-1)}
\int_0^1 \!\! \int_{\sph} t^{ \frac{n+\beta-4}{\beta}}
 X(t)^{\frac{4-2n}{n}}  |\nabla_{\omega}w_t|^2 dS \, dt 
\]
and on the other hand, recalling \eqref{alpha},
\bean
&& \hspace{-1.2cm}\int_0^1\!\! \int_{\sph} |\nabla_{\omega}w|^2 
\frac{1}{g'}(\alpha\alpha'' r+\alpha\alpha')\, dS \, dt \\
&=& \frac{n(n-4)^2}{8(n-1)\beta^2}\int_0^1 \int_{\sph} t^{\frac{n-\beta-4}{\beta}}
  X(t)^{\frac{4-2n}{n}} |\nabla_{\omega}w|^2  \, dS \, dt \\
&& -\frac{(n-2)(n-4)}{2(n-1)\beta} \int_0^1 \int_{\sph} t^{\frac{n-\beta-4}{\beta}}
  X(t)^{\frac{4-n}{n}} |\nabla_{\omega}w|^2  \, dS \, dt \\
  && -\frac{3(n^2-4)}{8n(n-1)} \int_0^1 \int_{\sph} t^{\frac{n-\beta-4}{\beta}}
  X(t)^{\frac{4}{n}} |\nabla_{\omega}w|^2  \, dS \, dt .
\eean
Combining the above we obtain (ii).
The proof of (iii) is much simpler and is omitted. \endproof

To proceed we define
\[
G^{\#}(t) =  G(t) - \Big( \frac{n}{2(n-1)}\Big)^3 A \; , \qquad t\in (0,1),
\]
where we recall that $A$ has been defined in Lemma \ref{lem1}.
\begin{lemma}
\la{wwh}
Let $v\in\cic (B_1\setminus\{0\})$
and let $w$ be defined by \eqref{5}. There holds
\begin{align*}
J_{\rm R}[v]  &= \Big( \frac{n}{2(n-1)}\Big)^3
\cA_{\rm R}[w] \\
& \quad +
 \int_0^1 \!\! \int_{\sph} t^{\frac{\beta+n-4}{\beta}}w_t^2 G^{\#}(t) dS \, dt + \int_0^1 \!\! \int_{\sph} t^{\frac{-\beta+n-4}{\beta}}w^2 H(t) dS \, dt
.
\end{align*}
\end{lemma}
\proof This is a direct consequence of Lemma \ref{lem:mark} (i).\endproof

\begin{lemma}
Let $n\geq 5$. If
\be
\beta\geq \beta_n :=  n\Big(  \frac{n^2-4n+8}{4n^4-24n^3+83n^2-120 n+52}
\Big)^{1/2}
\la{bn}
\ee
then the function $G^{\#}(t)$ is non-negative in $(0,1)$.
\end{lemma}
\proof We first note that
\bea
G^{\#}(t) &=&
\frac{ n(n^2-4n+8)}{4(n-1)} X(t)^{\frac{4-2n}{n}} \!
 -\frac{n^3(n^2-4n+8)}{16(n-1)^3\beta^2} 
 + \frac{5n(n-2)(3n-2)}{16(n-1)^3} X(t)^2  \nonumber \\
&=:& p_1 X(t)^{\frac{4-2n}{n}} +p_2 + p_3 X(t)^2
\la{zes}
\eea
Now, it easily follows from \eqref{zes} that
$G^{\#}(t)$ is monotone decreasing in $ (0,1]$. Hence its minimum 
equal to
\[
p_1+p_2+p_3 = \frac{ n(4n^4-24n^3+83n^2-120 n+52)}{16(n-1)^3}
-\frac{n^3(n^2-4n+8)}{16(n-1)^3\beta^2},
\]
which is non-negative if $\beta\geq\beta_n$.
\endproof

\begin{lemma}
Let $n\geq 5$ and $\beta\geq \beta_n$.
For any $w\in \cic(0,1)$ there holds
\[
 \int_0^1 t^{\frac{\beta +n-4}{\beta}} G^{\#}(t)  w_t^2 dt +
\int_0^1 t^{\frac{-\beta +n-4}{\beta}}H^{\#}(t) w^2 dt \geq 0
\]
where
\begin{align*}
H^{\#}(t) &=  -\, \frac{n(n-4)^2(n^2-4n+8)}{16(n-1)\beta^2} X^{\frac{4-2n}{n}}
+ \frac{(n-2)(n-4)(n^2-4n+8)}{4(n-1)\beta} X^{\frac{4-n}{n}} \\
& \quad + \frac{n^3(n-4)^2(n^2-4n+8)}{64(n-1)^3\beta^4}
+ \frac{3(n^2-4)(n^2-4n+8)}{16n(n-1)} X^{\frac{4}{n}}  \\
& \quad - \frac{n(n-2)(15n^3-104n^2+256n-152)}{32(n-1)^3\beta^2} X^2  \\
& \quad - \frac{5n(n-2)(n-4)(3n-2)}{16(n-1)^3\beta} X^3 
+ \frac{45(n-2)^2(3n-2)^2}{n(n-1)^3}X^4 .
\end{align*}
\la{lem:weighted}
\end{lemma}
\proof
Let $r_1,r_2$ be real numbers to be fixed later.
We have
\bean
 0 &\leq &\int_0^1 t^{\frac{\beta +n-4}{\beta}}
 G^{\#}(t)\Big( w_t +
 \frac{r_1+r_2X(t)}{t}w\Big)^2dt  \\
 &=&  \int_0^1 t^{\frac{\beta +n-4}{\beta}}G^{\#}(t)w_t^2dt + 
 \int_0^1 \Big\{ t^{\frac{-\beta +n-4}{\beta}}  G^{\#}(t)
 (r_1^2 +2r_1r_2X +r_2^2X^2)  \\
&& \quad\quad  -\Big( t^{\frac{n-4}{\beta}}G^{\#}(t)\big(r_1+r_2X(t)\big)   \Big)_t \Big\}w^2 dt
\eean
Substituting from \eqref{zes} and carrying out the computations we arrive at
\begin{align*}
0\leq &  \int_0^1 t^{\frac{\beta +n-4}{\beta}}G^{\#}(t)w_t^2dt \; + \\
& \int_0^1 t^{\frac{-\beta +n-4}{\beta}} \bigg\{
p_1r_1( r_1 -  \frac{n-4}{\beta})X^{\frac{4-2n}{n}}
+p_1( 2r_1r_2 -r_2 \frac{n-4}{\beta} +\frac{2n-4}{n}r_1)X^{\frac{4-n}{n}} \\
& \qquad + p_2r_1 (r_1 -\frac{n-4}{\beta}) +p_1r_2( r_2 
+\frac{n-4}{n})X^{\frac{4}{n}} + p_2r_2(2r_1 - \frac{n-4}{\beta})X \\
& \qquad +\big(  p_2r_2^2 -p_2r_2 +p_3r_1^2 -p_3r_1\frac{n-4}{\beta}    \big)X^2 
 + \big( 2p_3r_1r_2 -2p_3r_1 - p_3r_2\frac{n-4}{\beta}  \big)X^3 \\
& \qquad  + (p_3r_2^2 -3p_3r_2)X^4 \bigg\}w^2 dt \, .
\end{align*}
We now choose
\[
r_1 = \frac{n-4}{2\beta} \; , \qquad r_2 = - \frac{3(n-2)}{2n}.
\]
The choice for $r_1$ minimizes the coefficient of the leading term
in the last integral; 
the parameter $r_2$ is less important and the choice is made for
convenience. Substituting we obtain
\begin{align*}
 0\leq &  \int_0^1 t^{\frac{\beta +n-4}{\beta}}G^{\#}(t)w_t^2dt \; + \\
& \int_0^1 t^{\frac{-\beta +n-4}{\beta}} \bigg\{
 -\frac{n(n-4)^2(n^2-4n+8)}{16(n-1)\beta^2} X^{\frac{4-2n}{n}}
+ \frac{(n-2)(n-4)(n^2-4n+8)}{4(n-1)\beta} X^{\frac{4-n}{n}} \\
& \qquad + \frac{n^3(n-4)^2(n^2-4n+8)}{64(n-1)^3\beta^4}
+ \frac{3(n^2-4)(n^2-4n+8)}{16n(n-1)} X^{\frac{4}{n}}  \\
& \qquad - \frac{n(n-2)(15n^3-104n^2+256n-152)}{32(n-1)^3\beta^2} X^2 
 - \frac{5n(n-2)(n-4)(3n-2)}{16(n-1)^3\beta} X^3 \\
& \qquad  + \frac{45(n-2)^2(3n-2)^2}{n(n-1)^3}X^4 \bigg\}w^2 dt \, .
\end{align*}
which is the stated inequality.
\endproof

We next define the positive constants
\bea
&& \gamma_1 = \frac{n^6(n-4)^2}{256(n-1)^4}  , \qquad
\gamma_2 = \frac{3n^2(n-2)(5n-6)(n^2-4n+8)}{128(n-1)^4}  , \nonumber \\
&& 
\gamma_3 =\frac{9(n-2)(3n-2)(5n-6)(7n-6)}{256(n-1)^4} .
\la{gamma}
\eea

\begin{lemma}
Let $n\geq 5$ and $\beta\geq \beta_n$. Let  $v\in\cic (B_1\setminus\{0\})$ and let $w$ be defined by \eqref{5}.
We then have
\begin{align*}
& J_{\rm R}[v]  +  \int_0^{\infty}\!\! \int_{\sph}
v^2 r^{-1} \Big( \frac{\gamma_1}{\beta^4} 
X(r)^{\frac{2(n-2)}{n-1}} -\frac{\gamma_2}{\beta^2} X(r)^{\frac{3n-4}{n-1}} 
+\gamma_3 X(r)^4 \Big)dS \, dt \\
 & \hspace{2cm} \geq   \Big( \frac{n}{2(n-1)}\Big)^3
\cA_{\rm R}[w].
\end{align*}
\la{gol}
\end{lemma}
\proof From Lemmas \ref{wwh} and \ref{lem:weighted} we have
\[
J_{\rm R}[v]  \geq  \Big( \frac{n}{2(n-1)}\Big)^3 \cA_{\rm R}[w]  
   + 
  \int_0^1\!\! \int_{\sph} t^{\frac{n-\beta-4}{\beta}}w^2 \big( H(t) -H^{\#}(t) \big) dS \, dt .
\]
But we easily see that
\[
\frac{n}{2(n-1)}(H(t)-H^{\#}(t)) =
 -\frac{\gamma_1}{\beta^4} + \frac{\gamma_2}{\beta^2} X(t)^2 - \gamma_3 X(t)^4 ,
\]
hence
\bean
&& \hspace{-.8cm} J_{\rm R}[v] + \frac{2(n-1)}{n}
  \int_0^1 \!\! \int_{\sph}
 t^{\frac{n-\beta-4}{\beta}}w^2 \big(  \frac{\gamma_1}{\beta^4} 
- \frac{\gamma_2}{\beta^2} X(t)^2 + \gamma_3 X(t)^4 \big) dS \, dt \\
&& \geq \Big( \frac{n}{2(n-1)}\Big)^3 \cA_{\rm R}[w].
\eean
We now express the double integral above in terms of the function $v$
using once again \eqref{5}. 
We note that for any $\sigma\geq 0$ we have
\[
 \int_0^1 t^{\frac{n-\beta-4}{\beta}}w^2 X(t)^{\sigma}dt
= \frac{n}{2(n-1)}
  \int_0^1  r^{-1}v^2 X(r)^{\frac{\sigma n +4(n-2)}{2(n-1)}} dr  \, .
\]
Applying this for $\sigma=0,2,4$ we obtain the required inequality.
\endproof

\noindent
{\bf\em Proof of Theorem \ref{thm1}.} Let $u\in \cic( \Omega)$. Without
loss of generality we may assume that $\Omega =B_1$
and that $u\in\cic(B_1\setminus\{0\})$. 
Let $v=|x|^{\frac{n-4}{2}}u$. By the discussion following
Lemma \ref{lem1}, the required inequality is written
\[
\frac{J_{\rm R}[v]  + \frac{n^2(n-4)^2}{16}\Int_0^1 \!\! \int_{\sph} 
r^{-1}v^2X(r)^{\frac{2(n-2)}{n-1}}dS \, dr
+ J_{\rm NR}[v]}{  \Big( \Int_0^1  \!\! \Int_{\sph}
r^{-1}X(r)^{ \frac{2n-4}{n-4}}|v|^{\frac{2n}{n-4}}dS \, dr  \Big)^{\frac{n-4}{n}}}
\geq S_{2,n}.
\]
We make the choice
\[
\beta =\frac{n}{2(n-1)}.
\]

We shall prove the following two inequalities where $v$ and $w$ are related
by the change of variables \eqref{5}:
\bea
&& \hspace{-1cm}J_{\rm R}[v]  + \frac{n^2(n-4)^2}{16}\Int_0^1 \!\! \int_{\sph} 
r^{-1}v^2X(r)^{\frac{2(n-2)}{n-1}}dS \, dr \geq 
 \Big( \frac{n}{2(n-1)}\Big)^3
\cA_{\rm R}[w]  \la{in1} \\
&& \hspace{-1cm} J_{\rm NR}[v]  \geq 
 \Big( \frac{n}{2(n-1)}\Big)^3
\cA_{\rm NR}[w] .
\la{in2}
\eea
We claim that if these are proved then the result will follow. 
Indeed, by Lemma \ref{lem:mark} (iii)
the Sobolev terms are related by
\[
\int_0^1 \!\! \int_{\sph} r^{-1}X(r)^{ \frac{2n-4}{n-4}}
|v|^{\frac{2n}{n-4}}dS \, dr
= \frac{2(n-1)}{n}\int_0^1 \!\! \int_{\sph}
 |w|^{\frac{2n}{n-4}} t^{ \frac{n-\beta}{\beta}} dS \, dt \, .
\]
Hence, applying Lemma \ref{lem1} we shall obtain
\bean
&&\hspace{-2cm} \frac{  J_{\rm R}[v] + \frac{n^2(n-4)^2}{16}\Int_0^1 \!\! \int_{\sph} 
r^{-1}v^2X(r)^{\frac{2(n-2)}{n-1}}dS \, dr  +  J_{\rm NR}[v]  }{  \Big( \Int_0^1 \!\!
\Int_{\sph} r^{-1}X(r)^{ \frac{2n-4}{n-4}}|v|^{\frac{2n}{n-4}}dS \, dr 
\Big)^{\frac{n-4}{n}}} \\
&\geq & \Big( \frac{n}{2(n-1)}\Big)^{\frac{4(n-1)}{n}}
\frac{ \cA_{\rm R}[w] +  \cA_{\rm NR}[w]}{  
\Big( \Int_0^1 \!\! \Int_{\sph}  |w|^{\frac{2n}{n-4}} t^{ \frac{n-\beta}{\beta}} dS \, dt \Big)^{\frac{n-4}{n}}} \\
&\geq & \Big( \frac{n}{2(n-1)\beta}\Big)^{\frac{4(n-1)}{n}}S_{2,n} \\
&=& S_{2,n},
\eean
and the proof will be complete.

\vspace{0.2cm}

\noindent
{\em Proof of \eqref{in1}.}
For the specific choice of $\beta$ we have
\bean
&& \hspace{.5cm}\frac{\gamma_1}{\beta^4} 
X(r)^{\frac{2(n-2)}{n-1}} -\frac{\gamma_2}{\beta^2} X(r)^{\frac{3n-4}{n-1}} 
+\gamma_3 X(r)^4  \\
&=& \frac{\gamma_1}{\beta^4}X(r)^{\frac{2(n-2)}{n-1}}
\Big( 1 -\frac{\gamma_2}{\gamma_1}\beta^2 
X(r)^{\frac{n}{n-1}} +
\frac{\gamma_3}{\gamma_1}\beta^4 
X(r)^{\frac{2n}{n-1}}\Big) \\
&=& \frac{n^2(n-4)^2}{16}X(r)^{\frac{2(n-2)}{n-1}}
\Big( 1 -\frac{3(n-2)(5n-6)(n^2-4n+8)}{2n^2(n-1)^2(n-4)^2}
X(r)^{\frac{n}{n-1}}  \\
&& \quad + \frac{9(n-2)(3n-2)(5n-6)(7n-6)}{16n^2(n-1)^4(n-4)^2} 
X(r)^{\frac{2n}{n-1}}\Big)
\eean
The function
\[
y\mapsto  1 -\frac{3(n-2)(5n-6)(n^2-4n+8)}{2n^2(n-1)^2(n-4)^2}
y  + \frac{9(n-2)(3n-2)(5n-6)(7n-6)}{16n^2(n-1)^4(n-4)^2} y^2
\]
is convex and its values at the endpoints $y=0$ and $y=1$
do not exceed one. Noting that $n/(2n-2) >\beta_n$ the
result follows by Lemma \ref{gol}. 

\vspace{0.2cm}

\noindent
{\em Proof of \eqref{in2}.} We recall that the functional $\cA_{\rm NR}[w]$
has been defined in Lemma \ref{lem1} and the
functional $J_{\rm NR}[v]$ is expressed in terms of the function $w$
in Lemma \ref{lem:mark}. 

We observe that the coefficients of the terms involving
$(\Delta_{\omega} w)^2$ in the two sides of \eqref{in2} are equal. The same is
true for the coefficients of the terms involving
$|\nabla_{\omega} w_t|^2$.
Hence the result will follow if we establish that
\[
K(t) \geq \Big( \frac{n}{2(n-1)}\Big)^3 \cdot \frac{2(n-4)}{\beta^4}
= \frac{4(n-1)(n-4)}{n}.
\]
Indeed, the first two terms of $K(t)$ are enough for this, that is there holds
\[
(n-1)(n-4)X(t)^{\frac{8-4n}{n}} -\frac{(n-1)(n-4)^2}{n}
X(t)^{\frac{4-2n}{n}} - \frac{4(n-1)(n-4)}{n} \geq 0 
\]
for all $t\in (0,1)$.
This completes the proof of the Rellich-Sobolev inequality of Theorem \ref{thm1}.

The sharpness of the constant $S_{2,n}$ in the
Rellich-Sobolev inequality follows easily
by concentrating near a point
$x_0\in\partial\Omega$ with $|x_0|=D$.
$\hfill\Box$

\

\noindent
{\bf Acknowledgement.} We wish to thank the anonymous
referee for a very careful reading of the initial manuscript and 
for bringing
to our attention the relevance of article \cite{Be}.

\end{document}